\documentclass[11pt]{amsart}
\usepackage{fullpage} 
\usepackage{color} 
\usepackage{array} 
 \usepackage[normalem]{ulem} 
 \useunder{\uline}{\ul}{}
 \usepackage{algorithm} 
 \usepackage{algpseudocode}
\usepackage{graphicx}
\usepackage{hyperref}
 \usepackage{tikz} 
 \usetikzlibrary{arrows} 

\newtheorem{global-theorem}{Theorem}
\newtheorem{theorem}{Theorem}[section]

\newtheorem{proposition}[theorem]{Proposition}
\newtheorem{corollary}[theorem]{Corollary}

\theoremstyle{definition}
\newtheorem{definition}[theorem]{Definition}
\newtheorem{example}[theorem]{Example}
\newtheorem{remark}[theorem]{Remark}

\newcommand{\Cech}{\v{C}ech{} }

\begin{document}

\title{Grout: A 1-Dimensional Substitution Tiling Space Program}
\author{Scott Balchin \& Dan Rust}
\thanks{The University of Leicester, University Road, Leicester, LE1 7RH, United Kingdom}

\begin{abstract}
We introduce a GUI fronted program that can compute combinatorial properties and topological invariants of recognisable and primitive symbolic substitutions on finite alphabets and their associated tiling spaces. We introduce theory from the study of aperiodic $1$-dimensional tilings along with pseudocode highlighting the algorithms that we have implemented into the GUI. Grout is written using C\texttt{++} and its standard library.  
\end{abstract}
\maketitle


\section*{Supplementary Resources}
Grout is available to download for Windows and Mac OSX along with the supporting documentation at the following URL

\vspace{3mm}

\centerline{\small{\url{www2.le.ac.uk/departments/mathematics/extranet/staff-material/staff-profiles/scott-balchin}}}

\section{Introduction to Tiling Spaces}
	The study of aperiodic tilings of the plane has a rich history which emerged from the worlds of computer science and mathematical logic when Berger proved the undecidability of the domino problem in the 1960s \cite{B:undecidability}. It is now also a topic of general interest to the study of dynamical systems, topology, Diophantine approximation, ergodic theory, computer graphics, mathematical physics and even virology \cite{CH:codim-one-attractors, S:book, BS:beta-numerations, BG:book, L:comp-sci-book, BHZ:gap-labeling, T:virology-tilings}. Once it was discovered that sets of tiles exist which can only tile the plane aperiodically, a flurry of new discoveries quickly followed, culminating in the celebrated discovery of the famous \emph{Penrose tilings} \cite{P:pentaplexity} and in the discovery of quasicrystals \cite{S:nobel} for which Schectman was awarded the Nobel prize in Chemistry in 2011.

	The seminal paper of Anderson and Putnam \cite{AP} lead to a new algebraic topological approach to aperiodic tilings associated to a particular method of generating tilings of the plane known as the \emph{substitution method}. It is this method, but restricted to a 1-dimensional analogue, which we address here.
	
	It is hoped that the use of this program will make testing conjectures in tiling theory and symbolic substitutional dynamics more efficient, as well as allowing for the confirmation of hand calculations and comparing different methods of calculation (especially methods of calculating cohomology). Analysis of large data sets which can be potentially generated by the Grout source code, and the recognition of underlying patterns in the data may also aid to further the theory.
	
	The GUI front end for Grout is powered by Qt\cite{qt}. Grout has been designed with user experience in mind and includes many ease-of-use properties such as the ability to save and load examples, and convenient methods of sharing examples with other users via short strings that encode a substitution. There is also an option to export all of the data that has been calculated to a pre formatted \LaTeX{} file including all the TikZ code for the considered complexes.  This should be useful for those needing to typeset such diagrams in the future by fully automating the generation of diagrams in TikZ.
	
	In Section \ref{Sec:background} we introduce basic notions of tiling spaces associated to substitutions on finite alphabets. In Section \ref{Sec:main} we introduce the relevant tiling theory along with pseudocode for most of the non-trivial components of Grout that have been implemented. Section \ref{Sec:cohomology} will cover specifically those methods implemented to compute cohomology for tiling spaces. Throughout, we give instances of these methods being applied to a well-known example substitution.  In Section \ref{res} we give a wide range of other examples of calculations from Grout.
	
	\textbf{Acknowledgements.} The authors would like to thank Fabien Durand for his helpful advice in piecing together a robust recognisability algorithm, and Etienne Pillin for helpful comments with regards to the coding.  We would also like to thank Alex Clark, Greg Maloney and Jamie Walton for helpful suggestions during the development of Grout and for testing early versions of the program.

	\subsection{Background}\label{Sec:background}
		\begin{definition}
			Let $\mathcal{A}=\{a^1, \dots , a^l\}$ be a finite alphabet on $l$ symbols, and for all positive integers $n$ let $\mathcal{A}^n$ be the set of length $n$ words in $\mathcal{A}$.  Denote the union of these as $\mathcal{A}^+ = \bigcup_{n \geq 1} \mathcal{A}^n$. If, further, the empty word $\epsilon$ is included, we denote the union $\mathcal{A}^+ \cup \{\epsilon\}$ by $\mathcal{A}^\ast$.  A \emph{substitution} $\phi$ on $\mathcal{A}$ is a function $\phi\colon \mathcal{A} \to \mathcal{A}^+$ which assigns to each letter $a$ in $\mathcal{A}$ a non-empty word $\phi(a)$ whose letters are elements in $\mathcal{A}$. We extend $\phi$ to a function $\phi\colon \mathcal{A}^+ \to \mathcal{A}^+$ by concatenation; given a word $w=w_1 \dots w_n \in \mathcal{A}^n$, we set $\phi(w)=\phi(w_1)\dots \phi(w_n)$. In this way, we can consider finite iterates of the substitution $\phi^n$ acting on $\mathcal{A}^+$.
		\end{definition}
		
		\begin{definition}
			Let $\phi\colon \mathcal{A} \to \mathcal{A}^+$ be a substitution. We say a word $w \in \mathcal{A}^\ast$ is \emph{admitted} by the substitution $\phi$ if there exists a letter $a\in\mathcal{A}$ and a natural number $n\geq 0$ such that $w$ is a subword of $\phi^n(a)$ and denote by $\mathcal{L}^n\subset\mathcal{A}^n$ the set of all words of length $n$ which are admitted by $\phi$. Our convention is that the empty word $\epsilon$ is admitted by all substitutions. We form the \emph{language} of $\phi$ by taking the set of all admitted words $\mathcal{L} = \bigcup_{n\geq 0} \mathcal{L}^n$.
			
			We say a bi-infinite sequence $s \in \mathcal{A}^\mathbb{Z}$ is \emph{admitted} by $\phi$ if every subword of $s$ is admitted by $\phi$ and denote by $X_\phi$ the set of all bi-infinite sequences admitted by $\phi$. The symbol $s_i$ denotes the label assigned to the $i$th component of the sequence $s$. The set $X_\phi$ has a natural (metric) topology inherited from the product topology on $\mathcal{A}^\mathbb{Z}$ and a natural shift map $\sigma\colon X_\phi \to X_\phi$ given by $\sigma(s)_i= s_{i+1}$. We call the pair $(X_\phi, \sigma)$ the \emph{subshift} associated to $\phi$ and we will often abbreviate the pair to just $X_{\phi}$ when the context is clear.
		\end{definition}
			We say $\phi$ is a \emph{periodic} substitution if $X_\phi$ is finite, and say $\phi$ is aperiodic otherwise. If $\phi$ is \emph{aperiodic} and has a property which will be introduced later known as \emph{primitivity}, then $X_\phi$ is a Cantor set (in particular $X_\phi$ is non-empty) and $\sigma$ is a minimal action on $X_\phi$ - that is, the only closed shift-invariant subsets of $X_\phi$ are the empty set $\emptyset$ and the subshift itself $X_\phi$. Equivalently, the orbit of every point under $\sigma$ is dense in $X_\phi$.
		
		\begin{definition}
			Let $\phi$ be a substitution on the alphabet $\mathcal{A}$ with associated subshift $X_\phi$. The \emph{tiling space} associated to $\phi$ is the quotient space $$\Omega_\phi = (X_\phi\times [0,1]) /{\sim}$$ where $\sim$ is generated by the relation $(s,0)\sim (\sigma(s),1)$.
		\end{definition}
		One should imagine the point $(s,t)\in\Omega_\phi$ as being a partition or tiling of the real line $\mathbb{R}$ into labelled unit-length intervals (called tiles), where the labels are determined by the letters appearing in $s$ and the origin of $\mathbb{R}$ is shifted a distance $t$ to the right from the left of the tile labelled by the letter $s_0$. Two tilings $T$ and $T'$ in $\Omega_\phi$ are considered $\epsilon$-\emph{close} in this topology if, after a translate by a distance at most $\epsilon$, the tiles around the origin in $T'-\epsilon$ within a ball of radius $1/\epsilon$ lie exactly on top of the tiles around the origin in $T$ within a ball of the same radius and share the same labels.
		
		If $\phi$ is primitive and aperiodic, then $\Omega_\phi$ is a compact connected metric space which fibers over the circle with Cantor set fibers. The natural translation $T\mapsto T+t$ for $t\in\mathbb{R}$ equips $\Omega_\phi$ with a continuous $\mathbb{R}$ action which is minimal as long as $\phi$ is primitive. In this respect, tiling spaces are closely related to the more well-known spaces, the solenoids.
		
		Associated to any topological space $X$ is a collection of groups $\check{H}^*(X)$ called the \emph{\Cech cohomology} of $X$. We refer the reader to \cite{B:bott-tu} for an introduction to \Cech cohomology. \Cech cohomology is an important topological invariant of tiling spaces and it is of general interest to be able to calculate and study these groups. Grout implements three different methods for calculating the cohomology of tiling spaces associated to symbolic substitutions on finite alphabets.
		\begin{enumerate}
			\item The method of Barge-Diamond complexes as introduced in \cite{BD}
			\item The method of Anderson-Putnam complexes as introduced in \cite{AP}
			\item The method of forming an equivalent \emph{left proper} substitution as outlined in \cite{DHS:return-words}
		\end{enumerate}
		All three outputs are algebraically equivalent -- that is, they represent isomorphic groups -- but it is not always obvious that this is the case given the presentations. This disparity between presentations of results for the equivalent methods was one of the major motivating factors for developing Grout.  These cohomologies are extremely laborious to calculate by hand for large alphabets unless special criteria are met.

\section{Grout and its functions}\label{Sec:main}
\subsection{Substitution Structure}

We begin by outlining how we encode a substitution rule into Grout and how we implement the substitution rule.  In general, we have done most of the implementation by string manipulation methods.

We use a class \emph{sub} which has as its element a vector of strings. We always assume that our alphabet is ordered $a, b, c, \dots$. The first entry of a \emph{sub} class vector is $\phi(a)$, the second is $\phi(b)$ and so on. To validate the input we check that the number of unique characters appearing in all of the $\phi(x)$ is equal to the length of the alphabet, which is the length of the vector.  The GUI also employs the use of regular expressions to prevent illegal characters from being entered.


\begin{example}[Fibonacci Substitution]
The Fibonacci tiling is given by a substitution rule on the alphabet $\mathcal{A}=\{a,b\}$ and is defined as

$$
	\phi\colon
	\begin{cases}
		a \mapsto & b\\
		b \mapsto & ba
	\end{cases}
$$

The Fibonacci substitution is our main example used throughout the paper. See section \ref{res} for a selection of outputs for other common examples of substitutions.

\begin{figure}[H]
\includegraphics[scale=0.5]{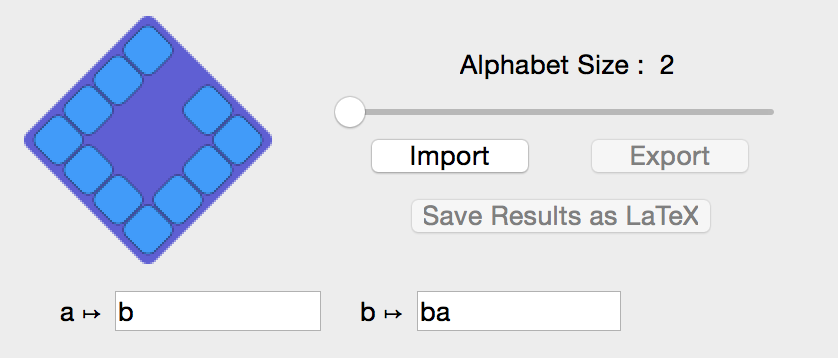}
\caption{The Fibonacci substitution entered into Grout}
\end{figure}

\end{example}

Next we implement a way to perform an iteration of $\phi$ on a string. We do not include any checks to validate that the string can be iterated on, as all strings that will be passed to this function will be created by the program itself, and therefore valid.

\begin{algorithm}
\caption{Substitution functions}
\label{alg2}
\begin{algorithmic}[1]
\\ \textbf{function} \underline{iterate}(string rhs)
\\ $\quad$ result = empty string
\\  $\quad$ \textbf{for} each character x in rhs \textbf{do}
\\     $\quad$ $\quad$       \textbf{append} $\phi$(x) to result
\\ $\quad$     \textbf{output} result
\end{algorithmic}
\end{algorithm}

%
%
%
%
%
%
%
%
%
%


\subsection{Substitution Matrices and their Properties}

\begin{definition}
Let $\mathcal{A} = \{a^1, \dots , a^l\}$ be an alphabet with a substitution $\phi \colon \mathcal{A} \to \mathcal{A}^+$, then $\phi$ has an associated \emph{substitution matrix} $M_\phi$ of dimension $l \times l$ given by setting $m_{ij}$ to be the number of times that the letter $a^i$ appears in $\phi(a^j)$.
\end{definition}

We will not give the algorithm of constructing the substitution matrix, the definition can be taken as a pseudo-algorithm.  We implement a square matrix class to work with the substitution matrix.  The first property that we will be checking for the substitution matrix is primitivity.

\begin{definition}
A substitution $\phi \colon \mathcal{A} \to \mathcal{A}^+$ is called \emph{primitive} if there exists a positive natural number $p$ such that the matrix $M^p_\phi$ has strictly positive entries. Such a matrix $M$ is also called \emph{primitive}.  This condition is equivalent to asking that there is a positive natural number $p$ such that for all $a,a' \in \mathcal{A}$ the letter $a'$ appears in the word $\phi^p(a)$.
\end{definition}

We will use the first definition to check the primitivity of our substitutions.  To check this condition on $M_\phi$, we count the number of zeros in the matrix and if the number of zeros reaches $0$ then we can conclude that the substitution is primitive. If not, then square the matrix and recount the number of zeros. If the number of zeros does not change then we can conclude that the substitution is not primitive.  This means that this check always halts.

\begin{algorithm}
\caption{Primitivity check}
\label{alg3}
\begin{algorithmic}[1]
\\ \textbf{function} \underline{primitive}
\\ $\quad$ matrix = substitution matrix of $\phi$
\\ $\quad$ \textbf{while} true \textbf{do} $\{$
\\ $\quad$ $\quad$ a = number of zeros in matrix
\\ $\quad$ $\quad$ matrix = matrix $\times$ matrix
\\ $\quad$ $\quad$ b = number of zeros in matrix
\\ $\quad$ $\quad$ \textbf{if} a=b \textbf{and} a!=0
\\ $\quad$ $\quad$ $\quad$ \textbf{output} false
\\ $\quad$ $\quad$ \textbf{if} b=0
\\ $\quad$ $\quad$ $\quad$ \textbf{output} true
\\ $ \quad$ $\}$
\end{algorithmic}
\end{algorithm}

We will be checking primitivity for all substitutions before we do calculations on them as if the substitution is not primitive many of the methods will not work, or will return false positive results.  Grout will always display whether a given substitution is primitive or not, it can also output the substitution matrix if asked to do so.


%
%
%

The next thing that we can do with the substitution matrix is give the tile frequencies and tile lengths of the substitution.  This requires us to compute the eigenvalues of the matrix.  We have implemented the QR method for computing the eigenvalues (for example see \cite{golub1996matrix}).  This gives us approximations to the real eigenvalues, and for the complex ones we simply give the conjugate pairs by their absolute values, and we give the results to two decimal places.  The eigenvalues of a substitution matrix may be printed out by ticking the eigenvalues box.

\begin{proposition}[Perron-Frobenius]
	Let $M$ be a primitive matrix.
	\begin{enumerate}
	 \item[i] There is a positive real number $\lambda_{PF}$, called the \emph{Perron-Frobenius eigenvalue}, such that $\lambda_{PF}$ is a simple eigenvalue of $M$ and any other eigenvalue $\lambda$ is such that $|\lambda|<\lambda_{PF}$.
	 \item[ii] There exist left and right eigenvectors, called the \emph{left and right Perron-Frobenius eigenvectors}, $\mathbf{l}_{PF}$ and $\mathbf{r}_{PF}$ associated to $\lambda_{PF}$ whose entries are all positive and which are unique up to scaling.
	\end{enumerate}

\end{proposition}
Given the above theorem, it is natural to ask what information is contained in the PF eigenvalue and eigenvectors of $M_\phi$ for a primitive substitution $\phi$.

If we were to assign a length to the tiles labelled by each letter, then we would hope for such a length assignment to behave well with the given substitution. The left PF eigenvector offers a natural choice of length assignments. If we assign to the letter $a^i$ the length $(\mathbf{l}_{PF})_i$, the $i$th component of the left PF eigenvector, then we can replace our combinatorial substitution by a geometric substitution. This geometric substitution expands the tile with label $a^i$ by a factor of $\lambda_{PF}$ and then partitions this new interval into tiles of lengths and labels given according to the combinatorial substitution. In order to give a unique output, Grout normalises the left PF eigenvector so that the smallest entry is $1$.

The information contained in the right PF eigenvector is also useful. The right PF eigenvector, once normalised so that the sum of the entries is $1$, gives the relative frequencies of each of the letters appearing in any particular bi-infinite sequence which is admitted by $\phi$. That is, if $|w|$ is the length of the word $w$, $|w|_i$ is the number of times the letter $a^i$ appears in the word $w$, and letting $s_{[-k,k]} = s_{-k} \ldots s_{-1} s_0 s_1 \ldots s_k$, then $$\lim_{k\to\infty} |s_{[-k,k]}|/|s_{[-k,k]}|_i = (\mathbf{r}_{PF})_i$$ for any $s\in X_\phi$.
\begin{figure}[H]
\includegraphics[scale=0.5]{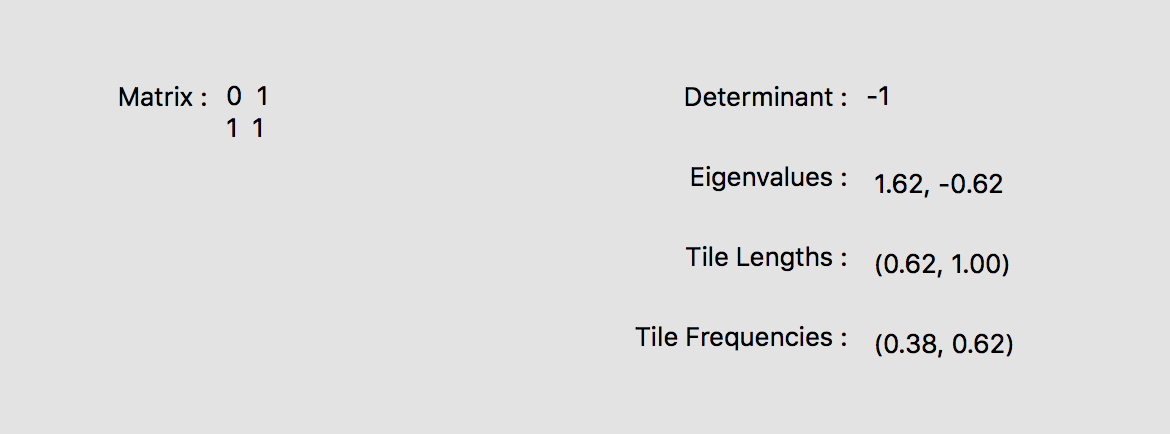}
\caption{The results for the matrix calculations for the Fibonacci substitutions}
\end{figure}

\subsection{Enumerating $n$-Letter Words}

Now that we have introduced the basic structure of the substitutions, and discussed the problem of primitivity and other matrix related calculations, we will discuss our first main function in Grout.

\begin{definition}
Given a substitution $\phi \colon \mathcal{A} \to \mathcal{A}^+$, we define the \emph{complexity function} at $n$ to be the number of unique $n$-letter words admitted by $\phi$.  We denote this function by $p_\phi(n)$, and so $p_\phi(n) = |\mathcal{L}_n|$.
\end{definition}

The complexity function of a tiling is a useful invariant \cite{J:complexity}. One is usually interested in either a deterministic formula for $p_\phi$ or information about the growth rate of $p_\phi$ such as polynomial degree; Grout can be used to at least give circumstantial evidence for these, though has no means of calculating either (this appears to be a very difficult problem in general).

Of particular interest are the number of 2 and 3 letter words, as we will be using them later to compute cohomology.  Our function will not only enumerate the number of $n$-letter words, but will also print out these words if required.  The algorithm uses C\texttt{++} sets as a data structure to store the $n$-letter words as it is automatically ordered and does now allow repetitions which leads to fast computation.  We start by generating a length $m$ admitted seed word $w$ such that $m\geq n$, and count all unique $n$-letter words appearing as subwords of the seed.  We then apply $\phi$ to the seed and add all new $n$-letter words to the result.  At each stage we count the size before and after adding the new words. If the size does not change we can stop, as no new $n$-letter words will be generated after a step without any new $n$-letter words.  It follows that the value $p_\phi(n)$ is computable in finite time for any fixed $n\geq 1$.

\begin{algorithm}
\caption{Finding all $n$-letter words}
\label{alg5}
\begin{algorithmic}[1]
\\ \textbf{function} \underline{nlw}(int n)
\\ $\quad$ result = empty ordered set
\\ $\quad$ seed = 'a'
\\ $\quad$ \textbf{while} seed length  $<$ n \textbf{do}
\\ $\quad$ $\quad$ seed = \underline{iterate}(seed)
\\
\\ $\quad$ difference = 1
\\ $\quad$ \textbf{while} difference != 0 \textbf{do} $\{$
\\ $\quad$ $\quad$ a = cardinality of result
\\ $\quad$ $\quad$ seed = \underline{iterate}(seed)
\\ $\quad$ $\quad$ \textbf{for} each $n$ length word $w$ in seed \textbf{do}
\\ $\quad$ $\quad$ $\quad$ \textbf{append} $w$ to result
\\ $\quad$ $\quad$ b = cardinality of result
\\ $\quad$ $\quad$ difference = b-a
\\ $\quad$ $\}$
\\ $\quad$ \textbf{output} result
\end{algorithmic}
\end{algorithm}

\begin{example}
It is well known that the complexity for the Fibonacci substitution satisfies $p_\phi(n)=n+1$, and we can verify this for any value of $n$ by computing its complexity in Grout.
%
%
%
%
%
%
%
%
%
%
\begin{figure}[H]
\includegraphics[scale=0.5]{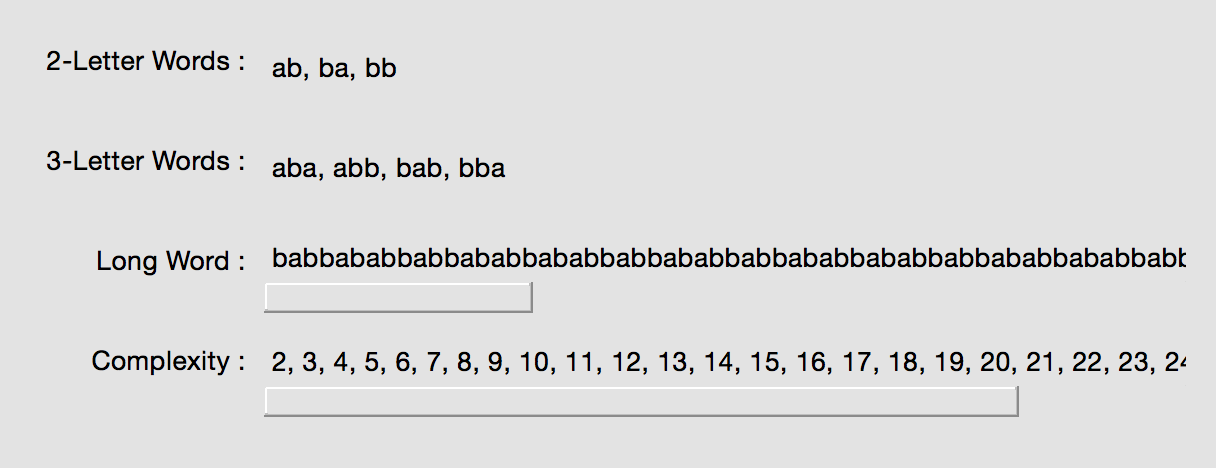}
\caption{The  words results for the Fibonacci substitution displayed in Grout}
\end{figure}
\end{example}

\subsection{Barge-Diamond and Anderson-Putnam Complexes}\label{sec:complexes}
Grout has the ability to output two simplicial complexes as PDFs (provided that the user has PDFLaTeX installed).  The first of these is the \emph{Barge-Diamond complex} \cite{BD}, which will be the key tool used in one of the methods of computing the \Cech cohomology of the tiling space.

\begin{definition}
Let $\mathcal{A} = \{a^1, \dots , a^l\}$ be an alphabet with a substitution $\phi \colon \mathcal{A} \to \mathcal{A}^+$, then we construct the \emph{Barge-Diamond complex} of $\phi$ as follows.  We have two vertices for each $a^i$, an \emph{in} node $v_i^+$ and an \emph{out} node $v_i^-$.  We draw an edge from $v_i^+$ to $v_i^-$ for all $i$.  Then for all two letter words $a^ia^j\in\mathcal{L}^2$ admitted by $\phi$, we draw an edge from $v_i^-$ to $v_j^+$.
\end{definition}

\begin{example}
As we have seen previously, the only two letter words admitted by the Fibonacci substitution are $ab,ba$ and $bb$, this gives us the following Barge-Diamond complex output in Grout.

\begin{figure}[H]
\centerline{
\includegraphics[scale=0.3]{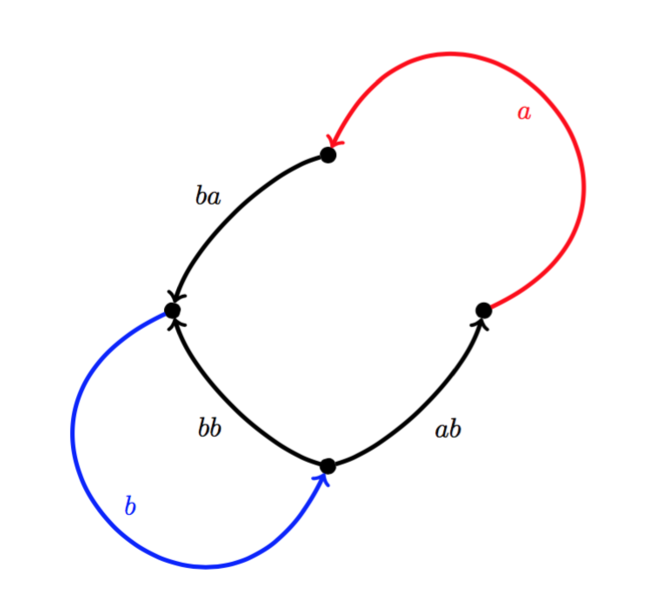}
}
\caption{The Barge-Diamond complex for the Fibonacci substitution.}
\end{figure}
\vspace{3mm}
\end{example}

The other complex that we consider for a substitution is (a variant of) the \emph{collared Anderson-Putnam complex} \cite{AP} which will again be the key tool used in one of the methods of computing \Cech cohomology of the tiling space. For brevity we will often shorten this to the \emph{AP complex}. This particular definition is based on what G\"{a}hler and Maloney call the \emph{Modified Anderson-Putnam complex} in \cite{GM:multi-one-d}.  The AP complex is constructed in a similar fashion to the Barge-Diamond complex, but makes use of both the two and three letter words.

\begin{definition}
Let $\mathcal{A} = \{a^1, \dots , a^l\}$ be an alphabet with a substitution $\phi\colon \mathcal{A} \to \mathcal{A}^+$, then we construct the \emph{Anderson-Putnam complex} of $\phi$ as follows.  We have a vertex $v_{ij}$ for each two letter word $a^ia^j\in\mathcal{L}^2$ admitted by $\phi$.  We draw an edge from $v_{ij}$ to $v_{jk}$ if and only if the three letter word $a^ia^ja^k$ is admitted by $\phi$.
\end{definition}

\begin{remark}
	One should note that this modified AP complex is slightly different to the definition originally introduced by Anderson and Putnam. In particular, the original definition distinguishes between different occurrences of a two letter word $a^ia^j$ if the occurrences of three letter words containing as a subword $a^ia^j$ do not overlap on some admitted four letter word. For example, if the language of a substitution included the two letter word $ab$, the three letter words $xab, yab, abw, abz$, and the four letter words $xabw, yabz$ but the words $xabz$ and $yabw$ did not belong to $\mathcal{L}$, then the original definition of the AP complex would have two instances of vertices with the label $ab$, say $(ab)_1$ and $(ab)_2$. In our definition, these vertices are identified, so that $(ab)_1 \sim (ab)_2 \sim ab$. An example of such a substitution is given by $\phi\colon a\mapsto bc,\: b\mapsto baab,\: c\mapsto caac$ where we label exactly one vertex with the label $aa$, but the original definition would require we include two distinct vertices labelled $(aa)_1$ and $(aa)_2$.
\end{remark}
In our discussion of cohomology calculated via AP complexes in Section \ref{sec:ap-cohomology}, we use this version of the AP complex to describe the performed calculations, and Grout implements this particular method. It would have been possible to use the original definition, or one of the many variant AP complexes that have been defined in the literature. There are at least three such variants discussed in \cite{GM:multi-one-d}, of varying complexities and situations in which they can be used.

\begin{example}\leavevmode

\begin{figure}[H]
\centerline{
\includegraphics[scale=1]{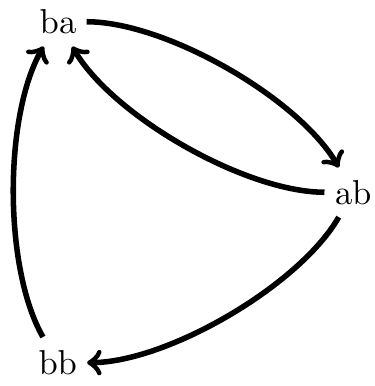}
}
\caption{The Modified Anderson-Putnam complex for the Fibonacci substitution}
\end{figure}

\vspace{3mm}
\end{example}

\subsection{Recognisability}
	\begin{definition}
		Let $\phi$ be a substitution on the alphabet $\mathcal{A}$. We say $\phi$ is \emph{recognisable} if for every bi-inifnite sequence $s\in X_\phi$ admitted by $\phi$ there is a unique way of decomposing $s$ into words of the form $\phi(a)$ for $a\in\mathcal{A}$. That is, up to a finite shift, there exists a unique bi-infinite sequence $\ldots a_{-1} a_0 a_1 \ldots$ such that $s = \ldots \phi(a_{-1}) \phi(a_0) \phi(a_1) \ldots$ and so we can \emph{recognise} which substituted letter each letter in $s$ has come from.
		
		Equivalently, we say $\phi$ is recognisable if there exists a natural number $K\geq 1$ such that for all admitted words $w\in\mathcal{L}$ with $|w|> 2K$, there exist unique words $x,y$ of length $|x|,|y|\leq K$ and a unique admitted word $u\in\mathcal{L}$ such that $w = x \phi(u) y$.
	\end{definition}
Recognisability is an important property of a tiling as many of the tools used to study the topology of the associated tiling spaces rely on recognisability as a hypothesis, much like primitivity. Recognisability of a primitive substitution is equivalent to aperiodicity of the subshift $X_\phi$ \cite{M:aperiodic}, and we make use of this result to decide recognisability.  The algorithm designed to determine if a given substitution is recognisable relies on finding a fixed letter and return words to that fixed letter.

\begin{definition}
Given a substitution $\phi$ on an alphabet $\mathcal{A}$, the letter $a$ is said to be \emph{fixed} (on the left) of \emph{order} $k$ if there exists some integer $k$ such that $\phi^k(a) = a u$ for some word $u$. Every substitution has at least one fixed letter and the value of $k$ for such a letter is bounded by the size of the alphabet.
\end{definition}

\begin{definition}
Given a fixed letter $a$, a \emph{return word} to $a$ is a word $v$ such that $v = au $ for some (possibly empty) word $u \in (\mathcal{A}\setminus\{a\})^\ast$, and $va$ is an admitted word of the substitution.
\end{definition}

If $\phi$ is primitive then, due to the primitivity of the substitution, the set of return words to any letter is finite. This is a consequence of the \emph{linear recurrence} of the subshift $X_\phi$, shown by Damanik and Lenz \cite{DL:linearly-recurrent}. We omit a definition of linear recurrence here.

\begin{algorithm}
\caption{Finding all return words to fixed letter f}
\label{alg6}
\begin{algorithmic}[1]
\\ \textbf{function}\underline{ returnwords}(character f)
\\ $\quad$ result = empty ordered set
\\ $\quad$ length = 2
\\ $\quad$ \textbf{while} new return words are being added \textbf{do} $\{$
\\ $\quad$ $\quad$ nwords = \underline{nlw}(length)
\\ $\quad$  $\quad$ \textbf{for} all words $w$ in nwords \textbf{do} $\{$
\\ $\quad$ $\quad$ $\quad$ \textbf{if} last character of $w$ = first character of $w$= f \textbf{and} $w$ has no other f appearing \textbf{do}
\\ $\quad$ $\quad$ $\quad$ $\quad$ \textbf{append} $w$ to result
\\ $\quad$ $\quad$ $\}$
\\ $\quad$ $\quad$ length = length + 1
\\ $\quad$ $\}$
\\ $\quad$ \textbf{output} result
\end{algorithmic}
\end{algorithm}

We will use these return words to determine whether a substitution is recognisable or not. The following proposition appears in \cite{HL:periodicity}.
	\begin{proposition}
		Let $\phi$ be a primitive substitution on $\mathcal{A}$ and let $a$ be a fixed letter. Let $\mathcal{R}$ be the set of all return words to $a$. So $\mathcal{R} = \{v\mid v = au,\: aua\in\mathcal{L},\: u \in (\mathcal{A}\setminus\{a\})^\ast\}$. The substitution $\phi$ is not recognisable if and only if, for all $v,v'\in\mathcal{R}$, there exists a $p\geq 1$ such that $\phi^p(vv') = \phi^p(v'v)$.
	\end{proposition}

	As $\mathcal{R}$ is finite, and together with the next proposition which appears in \cite{C:iterates-of-words} and \cite{ER:iterates-of-words}, this gives us a finite deterministic check for recognisability.
	\begin{proposition}
		Let $\phi$ be a substitution on $\mathcal{A}$ and let $|\mathcal{A}| = n$. For words $u,w\in\mathcal{A}^+$, there exists a $p\geq 1$ such that $\phi^p(u) = \phi^p(w)$ if and only if $\phi^n(u) = \phi^n(w)$.
	\end{proposition}
	
	That is, if some iterated substitution of $u$ and $v$ are ever equal, then their iterates must become equal at least by the $n$th iteration of the substitution, where $n$ is the size of the alphabet. In the algorithm, $k$ is taken to be the $k$ from the definition of the fixed letter, and $n$ is the size of the alphabet.

\begin{algorithm}
\caption{Recognisability check}
\label{alg7}
\begin{algorithmic}[1]
\\ \textbf{function} \underline{recognisable}
\\ $\quad$ rwords = \underline{returnwords}(f)
\\ $\quad$ \textbf{for} each word $w$ in rwords \textbf{do}
\\ $\quad$ $\quad$ \textbf{for} each word $v\neq w$ in rwords \textbf{do}
\\ $\quad$ $\quad$ $\quad$ \textbf{if}($\phi^{k \times n}$($w+v$) = $\phi^{k \times n}$($v+w$))
\\ $\quad$ $\quad$ $\quad$ $\quad$ \textbf{output} true
\\ $\quad$ \textbf{output} false
\end{algorithmic}
\end{algorithm}

\begin{figure}[H]
\includegraphics[scale=0.5]{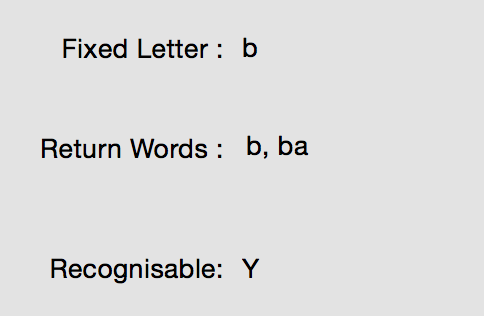}
\caption{The recognisability results for the Fibonacci substitution}
\end{figure}

\section{Cohomolgy of Tiling Spaces in Grout}\label{Sec:cohomology}
\subsection{Via Barge-Diamond}
	 Let $\phi$ be a primitive, recognisable substitution on the alphabet $\mathcal{A}$. Let $G$ be the Barge-Diamond complex of $\phi$ and let $S$ be the subcomplex of $G$ formed by all those edges labelled with two letter admitted words $a^ia^j$.
	
	Let $\tilde{\phi}\colon S \to S$ be a graph morphism defined in the following way on vertices. Let $l(i)$ and $r(i)$ be such that $\phi(a^i) = a^{l(i)} u a^{r(i)}$ and define $\tilde{\phi}(v_i^+) = v_{l(i)}^+$ and $\tilde{\phi}(v_i^-) = v_{r(i)}^-$. Note that if $a^ia^j$ is admitted by $\phi$, then $a^{r(i)}a^{l(j)}$ is also admitted by $\phi$, and so $\tilde{\phi}$ is a well defined graph morphism on $S$. As $\tilde{\phi}(S)\subset S$, we can define the eventual range $ER = \bigcap_{m\geq 0} \tilde{\phi}^m(S)$ (which stabilises after finitely many substitutions).
	
	For this method of computation we make use of the following result attributed to Barge and Diamond \cite{BD}.
	\begin{proposition}
		There is a short exact sequence
		$$0 \to \tilde{H}^0(ER) \to \varinjlim M_\phi^T \to \check{H}^1(\Omega_\phi) \to H^1(ER) \to 0.$$
	\end{proposition}
	The eventual range $ER$ is a (possibly disconnected) graph, and so $\tilde{H}^0(ER)$ and $H^1(ER)$ are finitely generated free abelian groups of some ranks $k$ and $l$ respectively. Hence we have
	\begin{corollary}\label{cor:BD}
		$$\check{H}^1(\Omega_\phi)\cong \varinjlim M_\phi^T/ \mathbb{Z}^k \oplus \mathbb{Z}^l.$$
	\end{corollary}
	
	Grout displays the \Cech cohomology using the Barge-Diamond method in the above form.
	
	These results fail in general if $\phi$ is not primitive or recognisable, and so Grout only performs this calculation after checking these two conditions.
	
	The only involved part of this calculation is finding the eventual range of the Barge-Diamond complex. After we have this we can simply find the number of connected components and use the Euler characteristic to find the reduced cohomology in rank zero and one.  Therefore, we now give the algorithm that we use to find the eventual range.  We denote by $w[0]$ the first letter of the word $ab$ and $w[1]$ the second letter of $ab$.
	
	\begin{algorithm}
\caption{Eventual range of the Barge-Diamond complex}
\label{alg8}
\begin{algorithmic}[1]
\\ \textbf{function} \underline{eventual range}
\\ $\quad$ twoletter = \underline{nlw}(2)
\\ $\quad$ difference = 1
\\ $\quad$ \textbf{while} difference != 0 \textbf{do} $\{$
\\ $\quad$ $\quad$  temp = empty ordered set
\\ $\quad$ $\quad$ \textbf{for} each word $w$ in twoletter \textbf{do}
\\ $\quad$ $\quad$ $\quad$ \textbf{append} last character of \underline{iterate}($w$[0]) + first character of \underline{iterate}($w$[1]) to temp
\\ $\quad$ $\quad$ difference = cardinality of twoletter - cardinality of temp
\\ $\quad$ $\quad$ twoletter = temp
\\ $\quad$ $\}$
\\ $\quad$ \textbf{output} twoletter
\end{algorithmic}
\end{algorithm}
	
\subsection{Via Anderson-Putnam}\label{sec:ap-cohomology}
	Let $\phi$ be a primitive, recognisable substitution on the alphabet $\mathcal{A}$. Let $K$ be the Anderson-Putnam complex of $\phi$. Anderson and Putnam showed in \cite{AP} that the \Cech cohomology of $\Omega_\phi$ is determined by the direct limit of an induced map acting on the cohomology of a CW complex. Using the modified AP complex $K$, G\"{a}hler and Maloney showed that this complex, as defined in Section \ref{sec:complexes}, can be used in place of the originally defined AP complex. We define the map acting on the AP complex $K$ in the following way.
	
	Again let $l(i)$ and $r(i)$ be such that $\phi(a^i) = a^{l(i)} u a^{r(i)}$. Let $E$ be an edge with label $a^ia^ja^k$ and define $L=|\phi(a^j)|$. Suppose $\phi(a^j) = a_1 a_2 \ldots a_L$. Define a continuous map called the \emph{collared substitution} $\tilde{\phi} \colon K \to K$ by mapping the edge $E$ to the ordered collection of edges with labels
	$$
	[a^{r(i)} a_1 a_2]
	[a_1 a_2 a_3]
	\cdots
	[a_{L-2} a_{L-1} a_L]
	[a_{L-1} a_{L} a^{l(k)}]
	$$
	in an orientation preserving way and at normalised speed. This map is well defined and continuous, hence an induced map $\tilde{\phi}^\ast \colon H^1(K) \to H^1(K)$ on cohomology exists. We use the following result from \cite{GM:multi-one-d} (Or \cite{AP} if $K$ is replaced with the original definition of the AP complex).
	\begin{proposition}
		$$\check{H}^1(\Omega_\phi) \cong \varinjlim (H^1(K), \tilde{\phi}^*)$$
	\end{proposition}
	Let $R=\mbox{rk} H^1(K)$, the rank of the cohomology of $K$. Grout finds an explicit generating set of cocycles for the cohomology of $K$ and then outputs the induced map $\tilde{\phi}^*$ as the associated $R\times R$-matrix $M_{AP}$, which should be interpreted as acting on $\mathbb{Z}^R$ with respect to this generating set. So $\check{H}^1(\Omega_\phi) \cong \varinjlim M_{AP}$.  
	
	The algorithm for this computation begins by constructing the boundary matrix for the AP-complex.  To do this we take all of the admitted three letter words $abc$ and we use the convention that the boundary of this edge is $bc-ab$.  Using this, we construct the associated $m \times n$ boundary matrix $B$ where $m$ is the number of two letter words and $n$ is the number of three letter words.  We then use standard methods from linear algebra to find a maximal set of linearly independent $n$ dimensional vectors $g$ such that $Bg = 0$, searching over the set of all vectors of $0$s and $1$s. By construction, this set generates the kernel of the boundary map inside the simplicial $1$-chain group of $K$. This gives us a generating set of cycle vectors for the first homology of $K$.
	
	We then apply the collared substitution to each of these generating vectors, giving us a new set of image vectors.  Using Gaussian elimination, we find the coordinates of these image vectors in terms of the generating vectors.  This induced map on homology can be represented as a square matrix. The transpose of this matrix $M_{AP}$ then represents the induced map on cohomology, and $M_{AP}$ is the output for the cohomology calculation via the Anderson-Putnam method.  It should be noted that this algorithm is not efficient in the case where the substitution has many three letter words, as the dimension $m$ of the $1$-chain complex is the dominant limiting factor when finding linearly independent generating cycles. The time complexity increases exponentially with respect to $m$.
	
\subsection{Via Properisation}
	For this method of computation we make use of a technique involving return words, as outlined in \cite{DHS:return-words}, for replacing a primitive substitution with an equivalent \emph{pre-left proper} primitive substitution. One may then use the fact that if $\phi$ is a recognisable pre-left proper primitive substitution, then $\check{H}^1(\Omega_\phi)\cong \varinjlim M_\phi^T$ (See \cite{S:book}).

	We begin by defining what it means to be proper.

	\begin{definition}\leavevmode
		\item A substitution is \emph{left proper} if there exists a letter $a\in\mathcal{A}$ such that the leftmost letter of $\phi(b)$ is $a$ for all $b \in \mathcal{A}$. That is, $\phi(b)=aw_b$ for some $w_b\in\mathcal{A}^+$.
		\item A substitution is \emph{right proper} if there exists a letter $a\in\mathcal{A}$ such that the rightmost letter of $\phi(b)$ is $a$ for all $b \in \mathcal{A}$. That is, $\phi(b)=w_ba$ for some $w_b\in\mathcal{A}^+$.
		\item A substitution is \emph{fully proper}\footnote{We will often abbreviate this to just \emph{proper}.} if it is both left and right proper.
		\item A substitution is \emph{pre-left proper} if some power of the substitution is left proper. Similarly for \emph{pre-right proper} and \emph{pre-fully proper}.
\end{definition}

	The following algorithm produces what we call the \emph{pre-left properisation} of a substitution, so-called because there exists a finite power of the new substitution which is left proper. As per usual $k$ will be the one from the definition of the fixed letter $f$. 

	Note that if $v$ is a return word to the fixed letter $a$, then $va\in\mathcal{L}$ and so $\phi^k(va)$ must also be admitted by $\phi$. But $\phi^k(va) = \phi^k(v) \phi^k(a)$, and both of $\phi^k(v)$ and $\phi^k(a)$ begin with the fixed letter $a$, hence $\phi^k(v)$ is an exact composition of return words to $a$. So, if we apply $\phi^k$ to a return word, then the result is a composition of return words.  We will denote by $\psi$ this newly constructed substitution rule on the new alphabet $\mathcal{R}$ of return words.

\begin{algorithm}
\caption{Pre-left properisation}
\label{alg9}
\begin{algorithmic}[1]
\\ \textbf{function} \underline{preprop}
\\ $\quad$ rwords = \underline{returnwords}(f)
\\ $\quad$ $\psi$ = empty substitution with alphabet size being the cardinality of rwords
\\ $\quad$ \textbf{for} all words $w$ in rwords \textbf{do} $\{$
\\ $\quad$ $\quad$ temp = $\phi^k$($w$)
\\ $\quad$ $\quad$ decomposition = decompose temp into return words $w_{i_1} \dots w_{i_m}$
\\ $\quad$ $\quad$ $\psi$($w$) = decomposition
\\ $\quad$ $\}$ 
\\ $\quad$ \textbf{output} $\psi$
\end{algorithmic}
\end{algorithm}

	This algorithm gives us a new substitution on possibly more letters than with what we began, and we may take a power of this substitution to get a left proper one. That such a power exists is clear. Indeed, every return word $v\in\mathcal{R}$ begins with the fixed letter $a$ and, by primitivity, $\phi^i(a)$ contains at least two copies of the letter $a$ for large enough $i$, so $\phi^i(a) = v_0 a u$ for some return word $v_0\in\mathcal{R}$ and some other word $u$. But then $\psi^i(v)$ must begin with $v_0$. It follows that $\psi^i$ is a left-proper substitution with leftmost letter $v_0$.

	We may also form an equivalent fully proper substitution on $\mathcal{R}$ by composing $\psi^i$ with its \emph{right conjugate}. The right conjugate $\phi^{(R)}$ of a left proper substitution $\phi$ is given by setting $\phi^{(R)}(b) = w_ba$ where $a$ is the fixed letter such that $\phi(b) = aw_b$ for all $b\in\mathcal{A}$. The right conjugate is a right proper substitution, and the composition of a left proper and right proper substitution is both left and right proper, hence fully proper.  It is easy to show (See \cite{DL:properisation}) that $X_{\phi\circ \phi^{(R)}}$ and $X_\phi$ are topologically conjugate subshifts. In fact, a word is admitted by $\phi\circ\phi^{(R)}$ if and only if it is admitted by $\phi$, so $X_{\phi\circ \phi^{(R)}}$ and $X_\phi$ are equal. Hence, $\check{H}^1(\Omega_{\phi\circ \phi^{(R)}}) \cong \check{H}^1(\Omega_\phi)$.
	
	We make use of the following which has been paraphrased from results appearing in the work of Durand, Host and Skau in \cite{DHS:return-words}.
	\begin{proposition}
		Let $\phi$ be a primitive substitution on $\mathcal{A}$ and let $\psi$ be the pre-left properisation of $\psi$. The tiling space $\Omega_\psi$ is homeomorphic to $\Omega_\phi$.
	\end{proposition}
	Hence we get the corollary
	\begin{corollary}
	$$\check{H}^1(\Omega_\psi) \cong \check{H}^1(\Omega_\phi)$$
	\end{corollary}
	As $\psi^i$ is left proper, $\check{H}^1(\Omega_{\psi^i})\cong \varinjlim M_{\psi^i}^T \cong \varinjlim (M_{\psi}^i)^T \cong \varinjlim M_{\psi}^T$. Hence $\check{H}^1(\Omega_\phi)\cong \varinjlim M_\psi^T$.

	Grout outputs the pre-left properisation $\psi$, the left properisation $\psi^i$, the full properisation $\psi^i \circ (\psi^i)^{(R)}$, and owing to the above, Grout also outputs the matrix $M_\psi^T$ in the cohomology section.

\begin{figure}[H]
\centerline{
\includegraphics[scale=0.5]{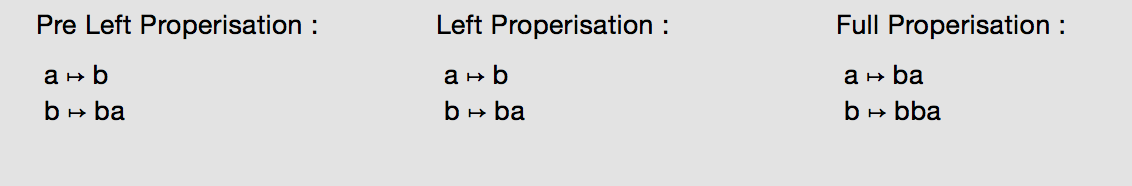}
}
\caption{The properisation results of the Fibonacci substitution}
\end{figure}

\begin{remark}
If the substitution is already proper, the properisation algorithm may return a different proper version of this substitution.  This may seem like it is a feature which has no use, but by iterating this process, we find that the sequence of substitutions is eventually periodic, first proved by Durand \cite{D:derived-seqs}. It was therefore decided to leave this feature intact, in order to study such sequences of properisations.
\end{remark}

\begin{figure}[H]
\centerline{
\includegraphics[scale=0.5]{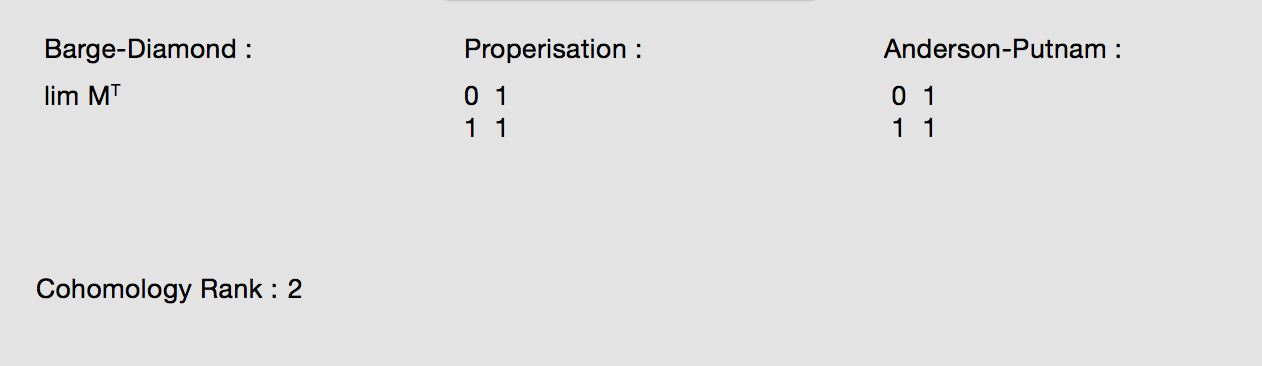}
}
\caption{The cohomology results of the Fibonacci substitution}
\end{figure}

\section{Further Results}\label{res}

In this section we provide the \LaTeX{}  output of cohomological calculations from Grout for a selection of both well-known and not so well-known substitutions appearing in the literature.

\begin{enumerate}
\item The \emph{Thue-Morse substitution} is one of the most well-studied substitutions in symbolic dynamics \cite{P:thue-morse}, possibly only superseded by the Fibonacci substitution in the attention it has received. It would be remiss to include a list of example outputs which did not include the results for this substitution.

\item The \emph{Tribonacci substitution} is an example of a unimodular irreducible Pisot substitution, first described by Rauzy in his seminal paper \cite{R:rauzy-fractals} introducing the so-called Rauzy fractals, and still actively studied for its interest to symbolic dynamicists and fractal geometers.

\item The \emph{Disconnected Subcomplex substitution} was first described by Barge and Diamond \cite{BD} as an example of a substitution whose Barge-diamond complex has a disconnected subcomplex of edges labelled by two letter words. This is reflected in the cohomology calculation via the Barge-Diamond complex, where a non-trivial quotient appears according to the formula for the cohomology described in Corollary \ref{cor:BD}.

\item The Fibonacci and Tribonacci substitutions are the first and second in an infinite family of primitive recognisable substitutions which take the form
$$
\begin{array}{rcl}
a_j & \mapsto & a_1a_{j+1},\:\:\: \mbox{ if } 1\leq j <n\\
a_n & \mapsto & a_1
\end{array}
$$
for an alphabet $\{a_1,\ldots a_n\}$ on $n$ letters. One might call these substitutions the \emph{$n$-ibonacci substitutions}. We have chosen to show the output for the Hexibonacci substitution where $n=6$, as the associated Barge-Diamond complex is particularly pleasing.
\end{enumerate}

\subsection{Thue-Morse}

\begin{center}

$$\begin{array}{rll}
 \color{red}a &\color{black}\mapsto& \color{red}a\color{blue}b\\
 \color{blue}b &\color{black}\mapsto& \color{blue}b\color{red}a
 \end{array}$$

\vspace{5mm}

\begin{tikzpicture}
\fill [red] (0.0,0.0) rectangle (0.5,0.5);
\fill [blue] (0.5,0.0) rectangle (1.0,0.5);
\fill [blue] (1.0,0.0) rectangle (1.5,0.5);
\fill [red] (1.5,0.0) rectangle (2.0,0.5);
\fill [blue] (2.0,0.0) rectangle (2.5,0.5);
\fill [red] (2.5,0.0) rectangle (3.0,0.5);
\fill [red] (3.0,0.0) rectangle (3.5,0.5);
\fill [blue] (3.5,0.0) rectangle (4.0,0.5);
\fill [blue] (4.0,0.0) rectangle (4.5,0.5);
\fill [red] (4.5,0.0) rectangle (5.0,0.5);
\fill [red] (5.0,0.0) rectangle (5.5,0.5);
\fill [blue] (5.5,0.0) rectangle (6.0,0.5);
\fill [red] (6.0,0.0) rectangle (6.5,0.5);
\fill [blue] (6.5,0.0) rectangle (7.0,0.5);
\fill [blue] (7.0,0.0) rectangle (7.5,0.5);
\fill [red] (7.5,0.0) rectangle (8.0,0.5);
\fill [blue] (8.0,0.0) rectangle (8.5,0.5);
\fill [red] (8.5,0.0) rectangle (9.0,0.5);
\fill [red] (9.0,0.0) rectangle (9.5,0.5);
\fill [blue] (9.5,0.0) rectangle (10.0,0.5);
\draw[ultra thick] (0,0) rectangle (10,0.5);
\end{tikzpicture}

Substitution Matrix : 
 
 \begin{equation*} \left( \begin{array}{cc} 
1 & 1\\1 & 1\\ \end{array} \right) 
 \end{equation*}

Full Properisation :
$$\begin{array}{rll}
a & \mapsto& cacbacab\\
b & \mapsto& cbabcacbabcbacab\\
c & \mapsto& cbabcacbacabcacbabcbacab
\end{array}$$
 
\vspace{5mm} 
 
Barge Diamond Cohomology Group : $\varinjlim M^T \oplus \mathbb{Z}^1$
 
 \vspace{2mm} 
 
Properisation Cohomology Matrix : 
 
 \begin{equation*} \left( \begin{array}{ccc} 
0 & 1 & 1\\1 & 0 & 1\\0 & 1 & 1\\ \end{array} \right) 
 \end{equation*} 
 
 \vspace{2mm}

Anderson-Putnman Cohomology Matrix : 
 
 \begin{equation*} \left( \begin{array}{ccc} 
0 & 0 & 0\\1 & 0 & 2\\1 & 1 & 1\\ \end{array} \right) 
 \end{equation*} 
 
 \vspace{2mm} 
 
Cohomology Rank : 2

\vspace{2mm}

\textbf{Barge-Diamond Complex}  
 
\vspace{-10mm} 
 
\begin{tikzpicture}[->, node distance=2cm, auto] 
\node [fill,circle,draw,inner sep = 0pt, outer sep = 0pt, minimum size=2mm] (ai) at (2.000000,0.000000) {}; 
\node [fill,circle,draw,inner sep = 0pt, outer sep = 0pt, minimum size=2mm] (ao) at (0.000593,2.000000) {}; 
\node [fill,circle,draw,inner sep = 0pt, outer sep = 0pt, minimum size=2mm] (bi) at (-2.000000,0.001185) {}; 
\node [fill,circle,draw,inner sep = 0pt, outer sep = 0pt, minimum size=2mm] (bo) at (-0.001778,-1.999999) {}; 
\draw (ai) edge[bend right=110, looseness=3, ->, red, ultra thick] node {$a$}(ao); 
\draw (bi) edge[bend right=110, looseness=3, ->, blue, ultra thick] node {$b$}(bo); 
\draw [-,ultra thick, bend right, draw=white, line width=6pt, looseness=0.7] (ao) to (ai); 
\draw [->,ultra thick, bend left, looseness=0.7] (ao) to node {$aa$} (ai); 
\draw [-,ultra thick, bend right, draw=white, line width=6pt, looseness=0.7] (bo) to (ai); 
\draw [->,ultra thick, bend right, looseness=0.7, swap] (bo) to node {$ab$} (ai); 
\draw [-,ultra thick, bend right, draw=white, line width=6pt, looseness=0.7] (ao) to (bi); 
\draw [->,ultra thick, bend right, looseness=0.7, swap] (ao) to node {$ba$} (bi); 
\draw [-,ultra thick, bend right, draw=white, line width=6pt, looseness=0.7] (bo) to (bi); 
\draw [->,ultra thick, bend left, looseness=0.7] (bo) to node {$bb$} (bi); 

 \end{tikzpicture}  
 
\end{center} 

\subsection{Tribonacci}

\begin{center}
  
$$\begin{array}{rll}
\color{red}a& \color{black}\mapsto& \color{red}a\color{blue}b\\
\color{blue}b& \color{black}\mapsto& \color{red}a\color{green}c\\
\color{green}c& \color{black}\mapsto& \color{red}a
\end{array}$$

\vspace{5mm}

\begin{tikzpicture}
\fill [red] (0.0,0.0) rectangle (0.5,0.5);
\fill [blue] (0.5,0.0) rectangle (1.0,0.5);
\fill [red] (1.0,0.0) rectangle (1.5,0.5);
\fill [green] (1.5,0.0) rectangle (2.0,0.5);
\fill [red] (2.0,0.0) rectangle (2.5,0.5);
\fill [blue] (2.5,0.0) rectangle (3.0,0.5);
\fill [red] (3.0,0.0) rectangle (3.5,0.5);
\fill [red] (3.5,0.0) rectangle (4.0,0.5);
\fill [blue] (4.0,0.0) rectangle (4.5,0.5);
\fill [red] (4.5,0.0) rectangle (5.0,0.5);
\fill [green] (5.0,0.0) rectangle (5.5,0.5);
\fill [red] (5.5,0.0) rectangle (6.0,0.5);
\fill [blue] (6.0,0.0) rectangle (6.5,0.5);
\fill [red] (6.5,0.0) rectangle (7.0,0.5);
\fill [blue] (7.0,0.0) rectangle (7.5,0.5);
\fill [red] (7.5,0.0) rectangle (8.0,0.5);
\fill [green] (8.0,0.0) rectangle (8.5,0.5);
\fill [red] (8.5,0.0) rectangle (9.0,0.5);
\fill [blue] (9.0,0.0) rectangle (9.5,0.5);
\fill [red] (9.5,0.0) rectangle (10.0,0.5);
\draw[ultra thick] (0,0) rectangle (10,0.5);
\end{tikzpicture}

Substitution Matrix : 
 
 \begin{equation*} \left( \begin{array}{ccc} 
1 & 1 & 1\\1 & 0 & 0\\0 & 1 & 0\\ \end{array} \right) 
 \end{equation*} 
 
\vspace{5mm}
 
Fixed Letter : a
 
 \vspace{2mm}Return Words : a, ab, ac
 
 \vspace{2mm}Recognisable: Yes 
 
 \vspace{2mm} 
 
Full Properisation : 
 
$$\begin{array}{rll}
a &\mapsto& bc\\
b &\mapsto& babc\\
c &\mapsto& bbc
\end{array}$$
 
\vspace{5mm} 
 
Barge Diamond Cohomology Group : $\varinjlim M^T$
 
 \vspace{2mm}

Properisation Cohomology Matrix : 
 
 \begin{equation*} \left( \begin{array}{ccc} 
0 & 0 & 1\\1 & 1 & 1\\0 & 1 & 0\\ \end{array} \right) 
 \end{equation*} 
 
 \vspace{2mm} 
 
Anderson-Putnman Cohomology Matrix : 
 
 \begin{equation*} \left( \begin{array}{ccc} 
0 & 0 & 1\\1 & 0 & 0\\1 & 1 & 1\\ \end{array} \right) 
 \end{equation*} 
 
 \vspace{2mm} 
 
Cohomology Rank : 3

\vspace{2mm}

\textbf{Barge-Diamond Complex}  

\vspace{-10mm}
 
\begin{tikzpicture}[->, node distance=2cm, auto] 
\node [fill,circle,draw,inner sep = 0pt, outer sep = 0pt, minimum size=2mm] (ai) at (2.000000,0.000000) {}; 
\node [fill,circle,draw,inner sep = 0pt, outer sep = 0pt, minimum size=2mm] (ao) at (1.000342,1.731853) {}; 
\node [fill,circle,draw,inner sep = 0pt, outer sep = 0pt, minimum size=2mm] (bi) at (-0.999316,1.732446) {}; 
\node [fill,circle,draw,inner sep = 0pt, outer sep = 0pt, minimum size=2mm] (bo) at (-2.000000,0.001185) {}; 
\node [fill,circle,draw,inner sep = 0pt, outer sep = 0pt, minimum size=2mm] (ci) at (-1.001368,-1.731260) {}; 
\node [fill,circle,draw,inner sep = 0pt, outer sep = 0pt, minimum size=2mm] (co) at (0.998289,-1.733038) {}; 
\draw (ai) edge[bend right=110, looseness=3, ->, red, ultra thick] node {$a$}(ao); 
\draw (bi) edge[bend right=110, looseness=3, ->, blue, ultra thick] node {$b$}(bo); 
\draw (ci) edge[bend right=110, looseness=3, ->, green, ultra thick] node {$c$}(co); 
\draw [-,ultra thick, bend right, draw=white, line width=6pt, looseness=0.7] (ao) to (ai); 
\draw [->,ultra thick, bend left, looseness=0.7] (ao) to node {$aa$} (ai); 
\draw [-,ultra thick, bend right, draw=white, line width=6pt, looseness=0.7] (bo) to (ai); 
\draw [->,ultra thick, bend right, looseness=0.7, swap] (bo) to node {$ab$} (ai); 
\draw [-,ultra thick, bend right, draw=white, line width=6pt, looseness=0.7] (co) to (ai); 
\draw [->,ultra thick, bend right, looseness=0.7, swap] (co) to node {$ac$} (ai); 
\draw [-,ultra thick, bend right, draw=white, line width=6pt, looseness=0.7] (ao) to (bi); 
\draw [->,ultra thick, bend right, looseness=0.7, swap] (ao) to node {$ba$} (bi); 
\draw [-,ultra thick, bend right, draw=white, line width=6pt, looseness=0.7] (ao) to (ci); 
\draw [->,ultra thick, bend right, looseness=0.7, swap] (ao) to node {$ca$} (ci); 

 \end{tikzpicture}

\end{center} 

\subsection{Disconnected Subcomplex}

\begin{center}
  
$$\begin{array}{rll}
\color{red}a& \color{black}\mapsto& \color{red}a\color{blue}b\color{green}c\color{magenta}d\color{red}a\\
\color{blue}b& \color{black}\mapsto& \color{red}a\color{blue}b\\
\color{green}c& \color{black}\mapsto& \color{green}c\color{magenta}d\color{blue}b\color{green}c\\
\color{magenta}d& \color{black}\mapsto& \color{magenta}d\color{blue}b
\end{array}$$

\vspace{5mm}

\begin{tikzpicture}
\fill [red] (0.0,0.0) rectangle (0.5,0.5);
\fill [blue] (0.5,0.0) rectangle (1.0,0.5);
\fill [green] (1.0,0.0) rectangle (1.5,0.5);
\fill [magenta] (1.5,0.0) rectangle (2.0,0.5);
\fill [red] (2.0,0.0) rectangle (2.5,0.5);
\fill [red] (2.5,0.0) rectangle (3.0,0.5);
\fill [blue] (3.0,0.0) rectangle (3.5,0.5);
\fill [green] (3.5,0.0) rectangle (4.0,0.5);
\fill [magenta] (4.0,0.0) rectangle (4.5,0.5);
\fill [blue] (4.5,0.0) rectangle (5.0,0.5);
\fill [green] (5.0,0.0) rectangle (5.5,0.5);
\fill [magenta] (5.5,0.0) rectangle (6.0,0.5);
\fill [blue] (6.0,0.0) rectangle (6.5,0.5);
\fill [red] (6.5,0.0) rectangle (7.0,0.5);
\fill [blue] (7.0,0.0) rectangle (7.5,0.5);
\fill [green] (7.5,0.0) rectangle (8.0,0.5);
\fill [magenta] (8.0,0.0) rectangle (8.5,0.5);
\fill [red] (8.5,0.0) rectangle (9.0,0.5);
\fill [red] (9.0,0.0) rectangle (9.5,0.5);
\fill [blue] (9.5,0.0) rectangle (10.0,0.5);
\draw[ultra thick] (0,0) rectangle (10,0.5);
\end{tikzpicture}

Substitution Matrix : 
 
 \begin{equation*} \left( \begin{array}{cccc} 
2 & 1 & 0 & 0\\1 & 1 & 1 & 1\\1 & 0 & 2 & 0\\1 & 0 & 1 & 1\\ \end{array} \right) 
 \end{equation*} 
 
 \vspace{5mm} 
 
Fixed Letter : a
 
 \vspace{2mm}Return Words : a, ab, abcd, abcdbcdb
 
 \vspace{2mm}Recognisable: Yes 
 
 \vspace{2mm}

Full Properisation : 
 
$$\begin{array}{rll}
a &\mapsto& cacad\\
b &\mapsto& cacabcad\\
c &\mapsto& cacaddbcad\\
d &\mapsto& cacaddbcaddbcabcad
\end{array}$$
 
\vspace{5mm} 
 
Barge Diamond Cohomology Group : $\varinjlim M^T / \mathbb{Z}^1$
 
 \vspace{2mm} 
 
Properisation Cohomology Matrix : 
 
 \begin{equation*} \left( \begin{array}{cccc} 
1 & 1 & 1 & 1\\0 & 1 & 0 & 1\\1 & 1 & 1 & 1\\0 & 0 & 1 & 2\\ \end{array} \right) 
 \end{equation*} 
 
 \vspace{2mm} 
 
Anderson-Putnman Cohomology Matrix : 
 
 \begin{equation*} \left( \begin{array}{cccc} 
1 & 0 & 1 & 0\\1 & 2 & 1 & 0\\0 & 1 & 1 & 1\\0 & 1 & 1 & 1\\ \end{array} \right) 
 \end{equation*} 
 
 \vspace{2mm} 
 
Cohomology Rank : 3
 
 \vspace{2mm}

\textbf{Barge-Diamond Complex} 

\vspace{-5mm}
 
\begin{tikzpicture}[->, node distance=2cm, auto] 
\node [fill,circle,draw,inner sep = 0pt, outer sep = 0pt, minimum size=2mm] (ai) at (2.000000,0.000000) {}; 
\node [fill,circle,draw,inner sep = 0pt, outer sep = 0pt, minimum size=2mm] (ao) at (1.414423,1.414004) {}; 
\node [fill,circle,draw,inner sep = 0pt, outer sep = 0pt, minimum size=2mm] (bi) at (0.000593,2.000000) {}; 
\node [fill,circle,draw,inner sep = 0pt, outer sep = 0pt, minimum size=2mm] (bo) at (-1.413585,1.414842) {}; 
\node [fill,circle,draw,inner sep = 0pt, outer sep = 0pt, minimum size=2mm] (ci) at (-2.000000,0.001185) {}; 
\node [fill,circle,draw,inner sep = 0pt, outer sep = 0pt, minimum size=2mm] (co) at (-1.415261,-1.413166) {}; 
\node [fill,circle,draw,inner sep = 0pt, outer sep = 0pt, minimum size=2mm] (di) at (-0.001778,-1.999999) {}; 
\node [fill,circle,draw,inner sep = 0pt, outer sep = 0pt, minimum size=2mm] (do) at (1.412746,-1.415680) {}; 
\draw (ai) edge[bend right=110, looseness=3, ->, red, ultra thick] node {$a$}(ao); 
\draw (bi) edge[bend right=110, looseness=3, ->, blue, ultra thick] node {$b$}(bo); 
\draw (ci) edge[bend right=110, looseness=3, ->, green, ultra thick] node {$c$}(co); 
\draw (di) edge[bend right=110, looseness=3, ->, magenta, ultra thick] node {$d$}(do); 
\draw [-,thick, bend right, draw=white, line width=4pt, looseness=0.7] (ao) to (ai); 
\draw [->,ultra thick, bend left, looseness=0.7] (ao) to (ai); 
\draw [-,thick, bend right, draw=white, line width=4pt, looseness=0.7] (bo) to (ai); 
\draw [->,ultra thick, bend right, looseness=0.7] (bo) to (ai); 
\draw [-,thick, bend right, draw=white, line width=4pt, looseness=0.7] (ao) to (bi); 
\draw [->,ultra thick, bend right, looseness=0.7] (ao) to (bi); 
\draw [-,thick, bend right, draw=white, line width=4pt, looseness=0.7] (co) to (bi); 
\draw [->,ultra thick, bend right, looseness=0.7] (co) to (bi); 
\draw [-,thick, bend right, draw=white, line width=4pt, looseness=0.7] (do) to (ci); 
\draw [->,ultra thick, bend right, looseness=0.7] (do) to (ci); 
\draw [-,thick, bend right, draw=white, line width=4pt, looseness=0.7] (ao) to (di); 
\draw [->,ultra thick, bend right, looseness=0.7] (ao) to (di); 
\draw [-,thick, bend right, draw=white, line width=4pt, looseness=0.7] (bo) to (di); 
\draw [->,ultra thick, bend right, looseness=0.7] (bo) to (di); 

 \end{tikzpicture} 
 
\end{center} 

\subsection{Hexibonacci}

\begin{center}

$$\begin{array}{rll}
\color{red}a& \color{black}\mapsto& \color{red}a\color{blue}b\\
\color{blue}b& \color{black}\mapsto& \color{red}a\color{green}c\\
\color{green}c& \color{black}\mapsto& \color{red}a\color{magenta}d\\
\color{magenta}d& \color{black}\mapsto& \color{red}a\color{brown}e\\
\color{brown}e& \color{black}\mapsto& \color{red}a\color{cyan}f\\
\color{cyan}f& \color{black}\mapsto& \color{red}a
\end{array}$$

\vspace{5mm}

\begin{tikzpicture}
\fill [red] (0.0,0.0) rectangle (0.5,0.5);
\fill [green] (0.5,0.0) rectangle (1.0,0.5);
\fill [red] (1.0,0.0) rectangle (1.5,0.5);
\fill [brown] (1.5,0.0) rectangle (2.0,0.5);
\fill [red] (2.0,0.0) rectangle (2.5,0.5);
\fill [blue] (2.5,0.0) rectangle (3.0,0.5);
\fill [red] (3.0,0.0) rectangle (3.5,0.5);
\fill [green] (3.5,0.0) rectangle (4.0,0.5);
\fill [red] (4.0,0.0) rectangle (4.5,0.5);
\fill [blue] (4.5,0.0) rectangle (5.0,0.5);
\fill [red] (5.0,0.0) rectangle (5.5,0.5);
\fill [magenta] (5.5,0.0) rectangle (6.0,0.5);
\fill [red] (6.0,0.0) rectangle (6.5,0.5);
\fill [blue] (6.5,0.0) rectangle (7.0,0.5);
\fill [red] (7.0,0.0) rectangle (7.5,0.5);
\fill [green] (7.5,0.0) rectangle (8.0,0.5);
\fill [red] (8.0,0.0) rectangle (8.5,0.5);
\fill [blue] (8.5,0.0) rectangle (9.0,0.5);
\fill [red] (9.0,0.0) rectangle (9.5,0.5);
\fill [cyan] (9.5,0.0) rectangle (10.0,0.5);
\draw[ultra thick] (0,0) rectangle (10,0.5);
\end{tikzpicture}

Substitution Matrix : 
 
 \begin{equation*} \left( \begin{array}{cccccc} 
1 & 1 & 1 & 1 & 1 & 1\\1 & 0 & 0 & 0 & 0 & 0\\0 & 1 & 0 & 0 & 0 & 0\\0 & 0 & 1 & 0 & 0 & 0\\0 & 0 & 0 & 1 & 0 & 0\\0 & 0 & 0 & 0 & 1 & 0\\ \end{array} \right) 
 \end{equation*} 
 
 \vspace{5mm} 
 
Fixed Letter : a
 
 \vspace{2mm}Return Words : a, ab, ac, ad, ae, af
 
 \vspace{2mm}Recognisable: Yes 
 
 \vspace{2mm}

Full Properisation : 
 
$$\begin{array}{rll}
a &\mapsto& bc\\
b &\mapsto& bdbc\\
c &\mapsto& bebc\\
d &\mapsto& bfbc\\
e &\mapsto& babc\\
f &\mapsto& bbc
\end{array}$$
 
\vspace{5mm} 
 
Barge Diamond Cohomology Group : $\varinjlim M^T$
 
 \vspace{2mm} 

Properisation Cohomology Matrix : 
 
 \begin{equation*} \left( \begin{array}{cccccc} 
0 & 0 & 0 & 0 & 0 & 1\\1 & 1 & 1 & 1 & 1 & 1\\0 & 1 & 0 & 0 & 0 & 0\\0 & 0 & 1 & 0 & 0 & 0\\0 & 0 & 0 & 1 & 0 & 0\\0 & 0 & 0 & 0 & 1 & 0\\ \end{array} \right) 
 \end{equation*} 
 
 \vspace{2mm}

Anderson-Putnman Cohomology Matrix : 
 
 \begin{equation*} \left( \begin{array}{cccccc} 
0 & 0 & 0 & 0 & 0 & 1\\1 & 0 & 0 & 0 & 0 & 0\\1 & 1 & 1 & 1 & 1 & 1\\0 & 0 & 1 & 0 & 0 & 0\\0 & 0 & 0 & 1 & 0 & 0\\0 & 0 & 0 & 0 & 1 & 0\\ \end{array} \right) 
 \end{equation*} 
 
 \vspace{2mm} 
 
Cohomology Rank : 6 
 
\vspace{5mm} 
 
\textbf{Barge-Diamond Complex} 

\vspace{-5mm}
 
\begin{tikzpicture}[->, node distance=2cm, auto] 
\node [fill,circle,draw,inner sep = 0pt, outer sep = 0pt, minimum size=2mm] (ai) at (2.000000,0.000000) {}; 
\node [fill,circle,draw,inner sep = 0pt, outer sep = 0pt, minimum size=2mm] (ao) at (1.732150,0.999829) {}; 
\node [fill,circle,draw,inner sep = 0pt, outer sep = 0pt, minimum size=2mm] (bi) at (1.000342,1.731853) {}; 
\node [fill,circle,draw,inner sep = 0pt, outer sep = 0pt, minimum size=2mm] (bo) at (0.000593,2.000000) {}; 
\node [fill,circle,draw,inner sep = 0pt, outer sep = 0pt, minimum size=2mm] (ci) at (-0.999316,1.732446) {}; 
\node [fill,circle,draw,inner sep = 0pt, outer sep = 0pt, minimum size=2mm] (co) at (-1.731557,1.000855) {}; 
\node [fill,circle,draw,inner sep = 0pt, outer sep = 0pt, minimum size=2mm] (di) at (-2.000000,0.001185) {}; 
\node [fill,circle,draw,inner sep = 0pt, outer sep = 0pt, minimum size=2mm] (do) at (-1.732742,-0.998802) {}; 
\node [fill,circle,draw,inner sep = 0pt, outer sep = 0pt, minimum size=2mm] (ei) at (-1.001368,-1.731260) {}; 
\node [fill,circle,draw,inner sep = 0pt, outer sep = 0pt, minimum size=2mm] (eo) at (-0.001778,-1.999999) {}; 
\node [fill,circle,draw,inner sep = 0pt, outer sep = 0pt, minimum size=2mm] (fi) at (0.998289,-1.733038) {}; 
\node [fill,circle,draw,inner sep = 0pt, outer sep = 0pt, minimum size=2mm] (fo) at (1.730963,-1.001881) {}; 
\draw (ai) edge[bend right=110, looseness=3, ->, red,ultra thick] node {$a$}(ao); 
\draw (bi) edge[bend right=110, looseness=3, ->, blue, ultra thick] node {$b$}(bo); 
\draw (ci) edge[bend right=110, looseness=3, ->, green, ultra thick] node {$c$}(co); 
\draw (di) edge[bend right=110, looseness=3, ->, magenta, ultra thick] node {$d$}(do); 
\draw (ei) edge[bend right=110, looseness=3, ->, brown, ultra thick] node {$e$}(eo); 
\draw (fi) edge[bend right=110, looseness=3, ->, cyan, ultra thick] node {$f$}(fo); 
\draw [-,thick, bend right, draw=white, line width=4pt, looseness=0.7] (ao) to (ai); 
\draw [->,ultra thick, bend left, looseness=0.7] (ao) to (ai); 
\draw [-,thick, bend right, draw=white, line width=4pt, looseness=0.7] (bo) to (ai); 
\draw [->,ultra thick, bend right, looseness=0.7] (bo) to (ai); 
\draw [-,thick, bend right, draw=white, line width=4pt, looseness=0.7] (co) to (ai); 
\draw [->,ultra thick, bend right, looseness=0.7] (co) to (ai); 
\draw [-,thick, bend right, draw=white, line width=4pt, looseness=0.7] (do) to (ai); 
\draw [->,ultra thick, bend right, looseness=0.7] (do) to (ai); 
\draw [-,thick, bend right, draw=white, line width=4pt, looseness=0.7] (eo) to (ai); 
\draw [->,ultra thick, bend right, looseness=0.7] (eo) to (ai); 
\draw [-,thick, bend right, draw=white, line width=4pt, looseness=0.7] (fo) to (ai); 
\draw [->,ultra thick, bend right, looseness=0.7] (fo) to (ai); 
\draw [-,thick, bend right, draw=white, line width=4pt, looseness=0.7] (ao) to (bi); 
\draw [->,ultra thick, bend right, looseness=0.7] (ao) to (bi); 
\draw [-,thick, bend right, draw=white, line width=4pt, looseness=0.7] (ao) to (ci); 
\draw [->,ultra thick, bend right, looseness=0.7] (ao) to (ci); 
\draw [-,thick, bend right, draw=white, line width=4pt, looseness=0.7] (ao) to (di); 
\draw [->,ultra thick, bend right, looseness=0.7] (ao) to (di); 
\draw [-,thick, bend right, draw=white, line width=4pt, looseness=0.7] (ao) to (ei); 
\draw [->,ultra thick, bend right, looseness=0.7] (ao) to (ei); 
\draw [-,thick, bend right, draw=white, line width=4pt, looseness=0.7] (ao) to (fi); 
\draw [->,ultra thick, bend right, looseness=0.7] (ao) to (fi); 

 \end{tikzpicture}  
 
\end{center}

\bibliographystyle{abbrv}
\bibliography{bib-grout}

\end{document}